\documentclass[12pt]{amsart}
\usepackage{amsmath,amscd,amsthm,amsfonts, amssymb,amsxtra}
\usepackage[all]{xy}
\newtheorem{thm}{Theorem}[section]  
\newtheorem*{un-no-thm}{Theorem}
\newtheorem{cor}[thm]{Corollary}     % Numbered along with thm
\newtheorem{lem}[thm]{Lemma}         % Numbered along with thm
\newtheorem{prop}[thm]{Proposition}  
\newtheorem{bigthm}{Theorem}
\newtheorem{bigcor}[bigthm]{Corollary}

\newtheorem*{mainthm}{Theorem}

\theoremstyle{definition}
\newtheorem{defn}[thm]{Definition}   % Numbered along with thm

\theoremstyle{definition}
\theoremstyle{definition}
\newtheorem*{ques}{Question}
\theoremstyle{remark}
\newtheorem{rem}[thm]{Remark}
\newtheorem{rems}[thm]{Remarks}

\newtheorem*{acks}{Acknowledgements}
 
\newtheorem*{out}{Outline}

\newtheorem{example}[thm]{Example}

\begin{document}
\title{The Moduli Space of Thickenings}
\date{\today}
\author{Mokhtar Aouina}
\address{Department of Mathematics, Jackson State University
Jackson, MS 39217}
\email{aouina@math.wayne.edu}

\begin{abstract} Fix $K$ a finite connected CW complex of dimension 
$\leq k$.  An {\it $n$-thickening} of $K$ is a pair $$(M,f)\, ,$$ 
in which $M$ is a compact $n$-dimensional manifold and $f\colon K \to M$ is a
simple homotopy equivalence. This concept was first introduced by
C.T.C.\ Wall approximately 40 years ago. Most of the known results
about thickenings are in a range of dimensions depending on $k$,
$n$ and the connectivity of $K$.

In this paper we remove the connectivity hypothesis on $K$. We
define {\it moduli space} of $n$-thickenings $T_n(K)$. We also
define a {\it suspension map} $E\colon T_n(K) \to T_{n{+}1}(K)$ and
compute its homotopy fibers in a range depending only on $n$ and
$k$. We will show that these homotopy fibers can be approximated
by certain section spaces whose definition depends only on the
choice of a certain stable vector bundle over $K$.
\end{abstract}

\maketitle%\pagestyle{fancy}
\def\Top{\bold T\bold o \bold p}
\def\wTop{\text{\rm w}\bold T}
\def\wT{\text{\rm w}\bold T}
\def\Sp{\bold S\bold p}
\def\vo{\varOmega}
\def\vs{\varSigma}
\def\smsh{\wedge}
\def\flush{\flushpar}
\def\id{\text{id}}
\def\dbslash{/\!\! /}
\def\codim{\text{\rm codim\,}}
\def\:{\colon}
\def\holim{\text{\rm holim\,}}
\def\hocolim{\text{\rm hocolim\,}}
\def\hodim{\text{\rm hodim\,}}
\def\hocodim{\text{hocodim\,}}
\def\Bbb{\mathbb}
\def\bold{\mathbf}
\def\Aut{\text{\rm Aut}}
\def\cal{\mathcal}
\def\Sec{\text{\rm sec}}
\def\Secst{\text{\rm sec}^{\text{\rm st}}}
\def\maps{\text{\rm map}}
%\setcounter{tocdepth}{1}
%\tableofcontents
%\addcontentsline{file}{sec_unit}{entry}
%\endtableofcontents

% Shortcut for "phi":
\newcommand{\f}{\varphi}
% Shortcut for a bar over a symbol ("top part"):
\newcommand{\tp}{\overline}
% Shortcut for derivation "delta":
\newcommand{\D}{d}
% Shortcut for category " Delta"
\newcommand{\Dc}{{\Delta}}
% Shortcut for standard n-simplicial or n-simplex and simplicial set " Delta"
\newcommand{\Ds}{\Delta}
%Shortcut for the category Sets
\newcommand{\Set}{\mathbf{Set}}

% For more than one paper by the same author in the bibliography:
%\providecommand{\bysame}{\makebox[3em]{\hrulefill}\thinspace}

% make obj hom op M (manifold) look vertical
\def\obj{\text{\rm obj\,}}
\def\hom{\text{\rm hom\,}}
\def\op{\text{\rm op\,}}
\def\M{\text{\rm M\,}}
\def\colim{\text{\rm colim\,}}
\def\maps{\text{\rm maps\,}}
\def\sec{\,\text{\rm  sec\,}}
\def\st{\text{\rm st\,}}
\def\:{\colon}
\def\cal{\mathcal}
\def\Bbb{\mathbb}
\def\bold{\mathbf}
\def\smsh{\wedge}
\def\Top{\bold T\bold o \bold p}
\def\Ho{\bold H\bold o}
\def\C_G{\Top^G}

\section{Introduction}\label{introduction}

 Let $K$ be a finite connected cell
complex of dimension $ \le k$.
An $\textit{n}$-{\it thickening}
of $K$ is a pair
$$(N,h)$$
in which $N$ is a compact smooth manifold and
$$
h\: K\to N\,
$$
is a simple homotopy equivalence. In addition, one assumes
\begin{itemize}
\item ({\it codimension $\ge 3$}). $k \le n{-}3$;
\item ({\it
$\pi$-$\pi$-condition}). The inclusion $\partial N \to N$ induces
an isomorphism of fundamental groups.
\item ({\it Surgery dimensions}). $n \ge 6$.
\end{itemize}

The above definition is due to C.T.C.\ Wall \cite{Wall1};
it generalizes the notion of a ``non-stable
neighborhood'' due to B.\ Mazur \cite{Mazur}.
\medskip

In any case, the $n$-thickenings of $K$ form the zero simplices
of a $\Ds$-set\footnote{i.e., simplicial set without the degeneracies.}
denoted by
$$
T_n(K)\, .
$$
A $1$-simplex of $T_n(K)$ can be thought of as an $s$-cobordism
between $n$-thickenings. More generally, the $j$-simplices of
$T_n(K)$ are ``$\Delta^j$-blocked'' families of thickenings. The
actual definition is too technical to give in this
introduction; we defer it to \S3.

The {\it moduli space} of $n$-thickenings $${\bold T}_n(K)$$ is
the geometric realization of $T_n(K)$. The main
purpose of this work is to obtain results about
the homotopy type of the moduli
space in a certain range of dimensions.
It will turn out that this range depends only on $n$ and $k$. In
\cite{Wall1}, Wall computed the set of path components of
the moduli space in a range that depends on $n$, $k$ and the {\it
connectivity} of $K$; our results enable one to dispense with
Wall's connectivity hypothesis.

We remark that the study of the path components of ${\bold T}_n(K)$
addresses a fundamental issue in differential topology:

\begin{ques}
{\it How does one enumerate the
set of compact $n$-manifolds
up to diffeomorphism within a fixed homotopy type?}
\end{ques}

In classical surgery theory, one imposes an
additional constraint:
the homotopy type of the boundary of the manifold
is held fixed. In the theory of thickenings,
this constraint is dropped.
\medskip

To study the thickening problem, Wall defined a {\it suspension map}
$$
E\:{\bold T}_n(K) \to {\bold T}_{n{+}1}(K)
$$
and studied the deviation from it being an isomorphism
on path components:

\begin{mainthm}[Wall's Suspension Theorem {\cite[\S 5]{Wall1}}]
 Assume
\begin{itemize}
\item $2n \geq
3k+3$, and
\item $K$ is a $(2k{-}n{+}1)$-connected based
CW complex of dimension $\le k$.
\end{itemize}
Then there is an exact sequence of sets
$$\{K\smsh K, S^n\}
\stackrel{P}\rightarrow \pi_0({\bold
T}_n(K))\stackrel{E}\rightarrow \pi_0({\bold
T}_{n+1}(K))\stackrel{H}\rightarrow \{K\smsh K, S^{n+1}\}\, .$$
\end{mainthm}

Here $\{X,Y\}$ denotes the abelian group of stable homotopy
classes of stable maps $X\to Y$. The names of the maps
 in the theorem are intended to remind one of the classical
homotopy theory EHP sequence (cf.\ \cite{James}).

In our initial study of ${\bold T}_n(K)$, we discovered a gap in
Wall's proof of the above theorem. The gap has to do with
exactness at the term $\pi_0({\bold T}_{n+1}(K))$. In
 a joint paper with J.R.\ Klein \cite{Aouina-Klein},
the gap was repaired.
The method we introduced in that paper ultimately led to
the generalization which appears in the present work.

We now wish to formulate our main results. This will
require some preparation.
Let $E\to B$ be a (Hurewicz) fibration
equipped with preferred section
$B\to E$. Denote the fiber over $b\in B$ by $F_b$.
Then one can form a new fibration
$$
Q_\bullet E \to B
$$
given by applying stable homotopy fiberwise. Explicitly, the fiber
$b\in B$ is $QF_b$, where $Q = \Omega^\infty \Sigma^\infty$
is the stable homotopy functor. We denote by
$$
\sec^{\rm st}(E\to B)
$$
the space of sections of $Q_\bullet E\to B$. This is the space of
{\it stable sections} of $E\to B$. It has the structure of an
infinite loop space; in particular, its set of path
components forms an abelian group.

We will apply this in the following situation: let $\xi$ be a rank
$n+1$ vector bundle over K, and let
$$
{\cal E}_\xi \to K \times K
$$
be the fibration whose fiber over $(x,y) \in K\times K$ is
given by the space
$$
S^\xi_x\smsh (\Omega_x^y K)_+
$$
where
\begin{itemize}
\item $S^\xi$ is the based $(n{+}1)$-sphere
given by taking the one point compactification
of the fiber of $\xi$ at $x\in K$;
\item $\Omega_x^y K$ is the space of paths from $x$ to $y$ in $K$;
\item $(\Omega_x^y K)_+$ is the effect of adding a basepoint to
$\Omega_x^y K $.
\end{itemize}

Let us say that an
$(n{+}1)$-thickening $\alpha = (N,h)$ {\it compresses} if there
exists an $n$-thickening $\beta=(M, g)$ such that $E(\beta)$ is in
the same path component as $\alpha$ in ${\bold T}_{n{+}1}(K)$.

To each such $\alpha = (N,h)$ we can associate a vector bundle
$\xi$ over $K$, which is given by pulling back the tangent
bundle of $N$ along the simple homotopy equivalence $h\: K \to N$.
We call $\xi$ the {\it tangential data} associated with $\alpha$.

Our first result determines the complete obstruction
to compressing an $(n{+}1)$-thickening in the metastable range.

\begin{bigthm}\label{compress}  To an $(n{+}1)$-thickening
 $\alpha = (N,f)$ with tangential data $\xi$,
one can assign an element
$$
e(\alpha) \in \pi_0(\sec^{\rm st}({\cal E}_\xi \to K \times K))
$$
called the {\bf Euler class}. This invariant vanishes when
$\alpha$ compresses.

Conversely, if $3k{+}1 \le 2n$, and the Euler class vanishes, then
$\alpha$ compresses.
\end{bigthm}

Our second and final main result concerns the determination
of the homotopy fibers of the suspension map $E$.

\begin{bigthm}\label{fiber_of_E}
 Assume $\alpha$ compresses.
Let ${\bold F}_n(\alpha)$ denote the homotopy fiber of $E\: {\bold T}_n(K) \to
{\bold T}_{n{+}1}(K)$ at $\alpha$. Then there is a
$(2n-3k-2)$-connected map
$$
{\cal H}\:{\bold F}_n(\alpha) \to \Omega \sec^{\rm st}({\cal E}_\xi \to
K\times K)\, .
$$
\end{bigthm}

With the assumptions of Theorem \ref{fiber_of_E},
the {\it Hopf invariant} is the
homotopy class represented by the composite
$$
H\: \Omega {\bold T}_{n{+}1}(K) \overset \partial \to {\bold
F}_n(\alpha) \overset {\cal H}\to
 \Omega \sec^{\rm st}({\cal E}_\xi \to K\times K)
$$
in which $\partial$ represents the boundary map in the long exact
sequence of the homotopy fibration, and $\Omega {\bold
T}_{n{+}1}(K)$ is the based loop space of ${\bold T}_{n{+}1}(K)$.

\begin{bigcor}
After fixing
a basepoint of  ${\bold T}_{n}(K)$ with tangential data $\xi$,
there is a long exact homotopy sequence {\tiny
$$
\xymatrix{ \pi_{2n{-}3k{-}1}(\sec^{\rm st}({\cal
E}_{\xi\oplus\epsilon} \to K\times K)) \ar[r]^-{P}
&\pi_{2n{-}3k{-}2}({\bold T}_n(K))
\ar[r]^E  & \pi_{2n{-}3k{-}2}({\bold T}_{n{+}1}(K)) \ar[dll]_-{H} \\
   \pi_{2n{-}3k{-}2}(\sec^{\rm st}({\cal E}_{\xi\oplus\epsilon} \to K\times K))\ar[r]^-{P} &
\quad \cdots \quad  \ar[r]^-{H} & \pi_1(\sec^{\rm st}({\cal
E}_{\xi\oplus\epsilon} \to K\times K))
 \ar[dl]_-{P} \\
 & \pi_0({\bold T}_{n}(K))  \ar[r]^E  &
\pi_0({\bold T}_{n{+}1}(K)) \, . }
$$}
\end{bigcor}

\begin{rems} (1). The last displayed
$H$ in the exact sequence isn't asserted to be homomorphism. This is
because the Hopf
invariant may fail to deloop.
\smallskip

{\flushleft (2).} The above is an exact
sequence of sets. In degrees $> 1$ it is an exact sequence
of abelian groups.

\smallskip
{\flushleft (3).} If in addition Wall's connectivity condition on
$K$ is assumed, then the above results reduce to Wall's suspension
theorem. The main
observation is that, in Wall's range, the  spaces
$$
\sec^{\rm st}({\cal E}_\xi \to K\times K) \quad \text{and} \quad
\text{\rm maps}^{\rm st}(K \smsh K,S^{n{+}1})
$$
have equivalent Postnikov sections in stage $1$ (where the right
side denotes the function space of stable maps from $K \smsh K$ to
$S^{n{+}1}$).
\end{rems}

\begin{out} \S2 is largely preliminary material on
$\Dc$-sets and manifold ads. In \S3 we define the moduli
space of thickenings. In \S4 we define the stabilization map.
\S5 and \S6 study the extent to which the stabilization map is
surjective on path components. In \S7 we indentify the homotopy
fibers of the stabilization map. In \S8 we compute some examples.
\end{out}

\begin{acks} This paper is a version, edited for publication,
of the author's Wayne State University  Ph.\ D.\ thesis. 
I am indebted to my advisor,
John R.\ Klein for choosing the topic. 
His help and insight was crucial in completing this work.
I am also indebted to Bob Bruner to introducing me to algebraic
topology and to Dan Isaksen for helping me identify mistakes.
\end{acks}

\section{Preliminaries} 

\subsection*{Spaces} We work in the category of compactly generated spaces.
A map $X\to Y$ of non-empty spaces is {\it $r$-connected} if its
homotopy fibers are $(r{-}1)$-connected
(by convention, every non-empty space is $({-}1)$-connected, the
empty space is $({-}2)$-connected, and a non-empty space $X$, is
$s$-connected for $s \ge 0$ if $\pi_j(X) = 0$ for $j \le s$ for
all choices of basepoint). An {\it weak (homotopy) equivalence}
is a map which is $\infty$-connected. When a map
$f\: X\to Y$ is a weak equivalence, we often indicate this
by $f\: \overset \sim \to Y$.

A {\it weak map} of spaces from $X$ to $Y$ consists of a finite chain
of maps
$$
X := X_0 \to X_1 \overset \sim \leftarrow X_2
\to
X_3 \overset \sim \leftarrow \quad \cdots \quad
\overset \sim \leftarrow X_{j{-}1} \to X_j =: Y
$$
in which each left pointing arrow is a weak equivalence.
A weak map induces an actual map in the homotopy category of spaces.
A weak map {\it $r$-connected} if all of the right pointing
maps are $r$-connected.
\bigskip

There are set theoretic difficulties
in defining the various moduli spaces appearing in this work.
The problem is related the fact that the
category of spaces is not small.
To get around this, we
fix a suitably large Grothendieck universe ${\cal U}$ and assume that
our spaces are subsets of that universe.
It will also be convenient to assume that
all manifolds are actually embedded in ${\Bbb R}^\infty$.

\subsection*{Review of $\Dc$-sets}
Introduced by Rourke and Sanderson \cite{Rourke-Sanderson}, 
$\Dc$-sets are simplicial sets without the
degeneracy operators. They arise naturally in classification
problems in differential topology, and enjoy many of the
properties satisfied by simplicial sets.

Let $\Dc$ be the category whose objects are {\it finite ordered
sets} and whose morphisms are the {\it order preserving
injections.} Note that each object of $\Dc$ with cardinality
$n{+}1$  is uniquely isomorphic to the ordered set
$$
[n] := \{0 < 1 <  \cdots < n\} \, .
$$
For this reason, a functor from $\Delta$ to any category is
determined by specifying a functor on the full subcategory of
$\Delta$ consisting  of the objects $[n]$.
%Let
%$$
%\widetilde{\Dc}
%$$
%be the category whose objects are the same as $\Dc$, and whose
%morphisms are the monotone functions (which we do not presume to
%be injections).

A {\it $\Delta$-set} is a functor
$$
X \colon \Dc^{\op}\rightarrow \Set \, ,
$$
where $\Set$ denotes the category of sets.

In $\Dc$, there are particular morphisms
$\delta_{i}\colon[n-1]\rightarrow[n]$, for $0\leq i\leq n$, given
by: \[
        \delta_{i}(j)=
        \begin{cases}
          j, &\text{if $j<i$;}\\
        j+1, &\text{otherwise.}\\
        \end{cases}
     \]

Let $X_{n}$ denote the value of this functor at the object $[n]$,
and $\D_i \colon X_n \to X_{n-1}$ be the map induced by
$\delta_i$. Then $X_{n}$ is called the set of {\it $n$-simplices},
and the function $\D_{i}$ is called the {\it $i^{th}$ face map}.
A 0-simplex is also called a vertex. The face maps satisfy the
identities
$$
\D_{i}\D_{j}=\D_{j-1}\D_{i}\,\,\,\,\,\text{if $\,\,i<j$.}
$$
Conversely, if we are given  sets $\{{X_{n}\}}_{n\in \mathbb{N}}$
and maps $\D_{i}$ satisfying the above identities, then
these define a unique $\Ds$-set.

Note that a simplicial set determines a $\Dc$-set by forgetting
its degeneracies.

%A {\it simplicial set}  is a functor
%$$
%\widetilde{\Dc}^{\op}\rightarrow \Set \, .
%$$
%The latter satisfies more identities corresponding to the fact
%that there are now degeneracies as well as face maps However, with
%$\Dc$-sets we work without degeneracies. Thus we have fewer
%relations to worry about.
\begin{example}
Let $\Ds[n]$ be the standard $n$-simplex. This is a
$\Ds$-set whose set of $m$-simplices $\Ds[n]_{m}$ is the set of
all monotone injective functions $f\colon [m]\rightarrow[n]$ and
whose face maps $\D_{i}\colon \Ds[n]_{m}\rightarrow \Ds[n]_{m-1} $
are given by $\D_{i}(f)=f\circ\delta_{i}$ for $0\leq i\leq m$.
\end{example}

\begin{defn}
\noindent
\begin{enumerate}
\item[1.] A {\it morphism} $f\:X \to Y$
 of $\Ds$-sets is just a natural transformation of functors.
%This is equivalent
%to the condition that $f$ satisfies the identities
%$f\partial_i = \partial_i
% f$.

\item[2.] A subfunctor $A \to X$ of $\Ds$-sets is called a {\it
sub $\Ds$-set}; this is sometimes written as a pair $(X,A)$.
\end{enumerate}

\end{defn}

\begin{example}
Let $\iota_{n}:=id_{[n]}$. The following are sub $\Ds$-sets of
$\Ds[n]$:

\begin{enumerate}

\item[1.] Let $\partial\Ds[n]$, the {\it boundary of $\Ds[n]$,}
be the $\Ds$-set whose
$j$-simplices are specified by:
 \[
        (\partial \Ds[n])_{j} =
        \begin{cases}
          \Ds[n]_{j}, &\text{if $0\leq j\leq n-1$;}\\
        \emptyset, &\text{otherwise.}\\
        \end{cases}
     \]

Thus $\partial\Ds[n]$ is the smallest sub $\Ds$-set of $\Ds[n]$
containing the faces $\delta_{j},\,0\leq j \leq n$ of the
$n$-simplex $\iota_{n}$.
%${\D\Ds[n]}_{j}=\Ds[n]_{j}\,\, \text{if $\,\,0\leq j\leq n-1$}$.

\item[2.] Let $\Lambda^k[n]$, called the {\it $k$-th  horn}, $\,0\leq k
\leq n$, be the sub $\Ds$-set of $\Ds[n]$, whose $j$-simplices are
given by:

\[
        (\Lambda^k[n])_{j}=
        \begin{cases}
          \Ds[n]_{j}, &\text{if $0\leq j\leq n-2$;}\\
          \Ds[n]_{n-1}-\{\delta_{k}\} &\text{if $j=n-1$;}\\
          \emptyset, &\text{otherwise.}\\
        \end{cases}
     \]

Thus $\Lambda^k[n]$ is a sub $\Ds$-set of $\Ds[n]$ which is
generated by all faces $\delta_{j}$ of the $n$-simplex $\iota_{n}$
except its $k^{th}$ face.
\end{enumerate}

\end{example}

We use $$\Dc \text{-} \Set$$ to denote the functor category of
$\Ds$-sets.

Let $\Delta^{n}$ be the standard affine topological $n$-simplex,
given by the convex hull of $[n]$ thought of as an orthonormal
basis of $\mathbb{R}^{n+1}$

A map $[m]\to [n]$ then induces an evident (linear) map of convex
hulls $a_*\:\Delta^m \to \Delta^n$.

\begin{defn}

The {\it (geometric) realization} of a $\Ds$-set $X$ is the
quotient space
$$
|X| = \frac{\coprod X_{n}\times \Delta^{n}}{\sim} ,
$$
in which $(a^*x,s)$ is identified with $(x,a_*t)$. Alternatively,
we can describe $|X|$  as the co-end of the functor $\Delta^\op
\times \Delta \to \Set$ given by $S \mapsto X_S \times \Delta^S$.

Note that $|X|$ has the structure of a CW complex
(\cite{Rourke-Sanderson}).

A map $X \to Y$ is a {\it weak equivalence} if it induces a weak
homotopy equivalence of spaces $|X| \to |Y|$.
\end{defn}

\subsubsection*{The extension condition}
Let $g\:X\to Y$ be a map of $\Ds$-sets. Consider a lifting problem
$$
\xymatrix{
 W \ar[r]^{\widehat{f}} \ar[d]_{\cap} & X  \ar[d]_g   \\
 Z \ar@.[ur] \ar[r]_f &Y\, .
}
$$
We say that $g$ is a {\it{Kan fibration}} if the lifting
problem has a solution for all such pairs of maps $f$ and
$\widehat{f}$ making the above diagram commute.

This is equivalent to asking that the lifting problem can be
solved when $(W,Z) = (\Delta^{n},\,\Lambda^i[n])$,
 for each $0\leq i \leq n$.

We say that $X$ is {\it Kan $\Ds$-set} if $X \to \ast$ is a Kan
fibration.

\subsubsection*{The Yoneda lemma}
The $j$-simplices of a $\Ds$-set $X$ are in bijection with the
$\Ds$-maps $\Ds[j]\to X$. The correspondence is given by
$$
(f\:\Ds[j]\to X) \,\, \mapsto \,\, f(\text{id}_{[j]}) \, ,
$$
where $\text{id}_{[j]}\: [j] \to [j]$ is the identity.

\subsubsection*{Homotopy groups}
Let $X$ be a $\Ds$-set equipped with basepoint $* \in X_0$.
The unique  map $a\:[j]\to [0]$ induces  a function $a\:X_0 \to X_j$
which in turn describes a basepoint of $X_j$. For this reason,
we can think of $X$ as a {\it based} $\Ds$-set.

Consider maps of $\Ds$-pairs.
$$
f\:(\Delta[j],\partial \Delta[j]) \to (X,*)
$$
Two such maps of $f$ and $g$ are said to be {\it homotopic} if
there is a map
$$
F\::\Delta[j]\times \Delta[1] \to X
$$
such that $F$ maps $\partial \Delta [j] \times \Delta[1]$ to the
basepoint, $F$ restricted to $\Delta[j]\times 0$ coincides with
$f$ and $F$ restricted to $\Delta[j] \times 1$ coincides with $g$.

If $X$ is a Kan $\Ds$-set, then homotopy is an equivalence relation.
In this case, we define
$$
\pi_j(X,*)
$$
to be the resulting set of equivalence classes. Then as
usual, $\pi_0(X,*)$ is a based set, $\pi_1(X,*)$ is group and
$\pi_j(X,*)$ is an abelian group when $j > 1$.

\subsubsection*{The function complex}
For $\Ds$-sets $X$ and $Y$ define
$$
F(X,Y)
$$
to be the $\Ds$-set whose $j$-simplices are
maps
$$
X\times \Ds[j]\to Y
$$
The face map $d_i\:F(X,Y)_j \to F(X,Y)_{j{-}1}$ is
given as follows: by the Yoneda lemma an element
of $F(X,Y)_j$ can be described as a certain kind of
function $X \to Y_j$ which we can post compose with
the $i$-th face map of $Y_j$
to get a function $X\to Y_{j-1}$. The latter then
corresponds to an element of $F(X,Y)_{j{-}1}$ (we omit the details).

The definition of the function complex generalizes to function spaces of
pairs, etc.\ in the evident way.

The function complex has the ``correct'' homotopy
type when $Y$ is Kan in the sense that the evident map
$$
|F(X,Y)| \to \text{maps}(|X|,|Y|)
$$
is a weak equivalence.

\subsubsection*{The homotopy fiber}
We will only need to describe the homotopy fiber
of an inclusion.

Let $(X,A)$ be a $\Ds$-pair,
 and assume that
$X$ is based. The {\it homotopy fiber} $F$  of $A\to X$ is the
$\Ds$-set whose $j$-simplices are given by the function complex of
maps of triads
$$
(\Ds[j] \times \Ds[1],\Ds[j] \times 1,\Ds[j]\times 0) \to (X,A,*)
$$
If $X$ is Kan, then $F$ has the correct homotopy type in the sense
that $|F|$ is identified with the homotopy fiber
at $*$ of the map $|A| \to |X|$.

\subsection*{(Very Special) Manifold Ads}
Let $n\ge 2$ be an integer. Let $M^m$ be a (smooth) manifold
equipped with submanifolds
$$
M(0),M(1),\dots,M(n{-}1) \,\, \subset \,\, \partial M
$$
which are codimension zero inside $\partial M$.
Set $[n-1]= \{0,1,...,n{-}1\}$.

We will assume:
\begin{itemize}
\item The total intersection of the $M_i$ is empty, and \item If
$i \ne j$, then $(M(i),\partial M(i))$ is transverse to
$(M(j),\partial M(j))$ and $M(i) \cap M(j)$ is a codimension zero
submanifold of $\partial M(j)$. \item For each nontrivial subset
$S \subset [n{-}1]$, let $M(S)$ denote the intersection of the
$M(j)$ for $j \in S$. Then $M(S)$ is a manifold and if $i \notin
S$, then $M(i)$ transversally intersects $M(S)$ in a codimension
zero submanifold of $\partial M(S)$.
\end{itemize}
Here  is some more notation: for $S \subset [n-1]$ a non-empty
subset, we define
$$
M_S := M([n{-}1]-S)\, .
$$
Our convention will be to set  $M_{[n-1]} := M$.

The above definition
is a special case of what is called an {\it $m$-manifold
$(n{+}1)$-ad.} This terminology, introduced by Wall
\cite[\S0]{Wall_book}, arises from the fact that one can
specify the data as an $(n{+}1)$-tuple
$$
(M;M(0),M(1),...,M(n-1))\, .
$$

\begin{example} A $3$-ad in the above sense consists of
$$
(M;M_1,M_0)
$$
in which $M_1 \amalg M_0$ is a
codimension zero compact submanifold of $\partial M$.
Define
$$
\partial_1 M \, \, := \,\, \text{closure}(\partial M - (M_0 \amalg M_1))
$$
Then $(M,\partial_1 M)$ is a cobordism between $(M_0,\partial M_0)$ and
$(M_1,\partial M_1)$.
\end{example}

We can also specify ads
in the language of functors.
Let ${\bold 2}^n$ be the poset of subsets of
$[n-1]$. Let ${\cal M}$ denote
be the category whose objects are smooth compact
manifolds with boundary and whose morphisms are
inclusions of one smooth manifold in the
boundary of another. Then above determines
an intersection preserving functor
$$
M_\bullet \: {\mathbf 2}^n \to {\mathcal M}
$$
defined by
$$
S\mapsto M_S \, .
$$

Here is an example of such an $(n{+}1)$-ad:
let $N$ be a compact manifold, possibly with
boundary. For each non-empty $S \subset  [n-1]$, form
$$
M_S\, \, := \,\, N \times {\Delta}^S\, .
$$
(We use the
standard method of ``rounding'' corners of $\Delta^{n-1}$
to consider $M_S$ as a smooth manifold with boundary).

%Let $M := M_S$ be a (smooth or PL) manifold equipped with
%codimension zero compact submanifolds $M_T \subset \partial M$,
%indexed by those subsets $T \subset S$ defined by deleting a
%single element of $S$. We will always assume the collection
%to be ``special'' in the sense that
%\begin{itemize}
%\item The intersection of the $M_T$, indexed over all
%$T \subset S$ given by deleting a single element,
%is empty.
%end{itemize}

%For each $U \subset T$, set
%$$
%M_U = \bigcap_{T \supset U} M_T \, .
%$$
%We will also require the following condition:
%\begin{itemize}
%\item for each pair of distinct subsets
% $U, V$ of $S$, the submanifolds
% $M_U$ and $M_V$ meet transversly. The boundaries of these
%manifolds are also required to meet transversely.
%\end{itemize}

For an $(n{+}1)$-ad $M_\bullet$, Set
$$
\partial_0 M_S := \coprod_{T \subsetneq S} M_T \, ,
$$
and set
$$
\partial_1 M_S := \text{closure}(\partial M_S - \partial_0 M_S) .
$$
Then we have a decomposition
$$
\partial M_S = \partial_0 M_S \cup \partial_1 M_S\, .
$$

\section{Thickenings}

Let $K$ be a finite connected cell complex of dimension $\le k$. An
{\it $n$-thickening} of $K$ is a pair
$$
(M,f)
$$
in which $M$ is an $n$-dimensional compact manifold with boundary
$\partial M$, and $f\:K \to M$ is a simple homotopy equivalence.
We additionally assume
\begin{itemize}
\item $k \le n{-}3$; \item $M$ satisfies the {\it
$\pi$-$\pi$-condition}, i.e., $\partial M$ is connected and the
inclusion $\partial M \to M$ induces an isomorphism of fundamental
groups.
\end{itemize}
A few words are in order. The first condition avoids knotting
phenomena. The second condition ensures that the map $\partial M
\to M$ is $(n{-}k{-}1)$-connected, which is a homotopy theoretic
way of thinking of the spine of $M$ as being $K$.

As promised in the introduction, we will now describe the
{\it moduli space} of $n$-thickenings of $K$. This will be the
geometric realization of a $\Dc$-set whose $j$-simplices are a
``$\Delta^j$-blocked'' version of a thickening.

\begin{defn}
Let $K^k$ be a connected finite cell complex of dimension $k$. We
denote by
$$
T_n(K)
$$
the $\Ds$-set in which a $0$-simplex is given by an
$n$-thickening $(M,f)$ of $K$.

If $j \ge 1$, then
$j$-simplex is specified by a pair
$$
(M_\bullet ,f_\bullet)
$$
in which $M_\bullet$ is  $(n+j)$-dimensional manifold $(j+2)$-ad
and
$$
f_\bullet \colon K\times \Delta^j\rightarrow M_\bullet
$$
is a simple $j$-blocked homotopy equivalence. This means that
$f_\bullet$ respects faces and for each nonempty subset $S \subset
[j]$, the map
 $$
f_S \colon K\times \Delta^S\rightarrow M_S
$$
is a simple homotopy equivalence. We also require
\begin{itemize}
\item the inclusion
$$
\partial_1 M_S \to M_S
$$
is an isomorphism of path components and also of fundamental
groups. \item For each $V \subset U$, the map
$$(M_{V}, \partial_1 M_{V})\subset (M_{U},\partial_1 M_{U})$$
is a simple homotopy equivalence of pairs.
\end{itemize}
\end{defn}

\begin{example} A $1$-simplex of $T_n(K)$
amounts to an $s$-cobordism $(M_{01},\partial_1 M_{01})$ between
$(M_0,\partial M_0)$ and $(M_1,\partial M_1)$, together with a
simple homotopy equivalence
$$
(K \times [0,1], K \times 0,K \times 1) \overset \sim \to
(M_{01},M_0,M_1)\, .
$$
\end{example}

\begin{defn}
Let $g\:X\to Y$ be a map of $\Ds$-sets. Consider a lifting problem

\begin{equation}
%\SelectTips{cm}{}
\xymatrix{
 W \ar[r]^{\widehat{f}} \ar[d]_{\cap} & X  \ar[d]^g   \\
 Z \ar@.[ur] \ar[r]_f &Y\, .
}
\end{equation}
We say that $g$ is a {\it{Kan fibration}} if the lifting
problem has a solution for all such pairs of maps $f$ and
$\widehat{f}$ making the above diagram commute.

This is equivalent to asking that the lifting problem can be
solved when $(W,Z) = (\Delta^{n},\,\Lambda^i[n])$,
 for each $0\leq i \leq n$.

We say that $X$ is {\it Kan $\Ds$-set} if $X \to \ast$ is a Kan
fibration.
\end{defn}

\begin{prop}
${T_{n}(K)}$ is a Kan $\Ds$-set.
\end{prop}

\begin{proof} A map $\Lambda^k[j] \to T_n(K)$ amounts
to specifying $j$ simplices of $T_n(K)_{j-1}:$
$$f_i\colon K\times
\Delta_i^{j-1}\rightarrow
M_i\,\,\,\,i=0,\,\cdots,\,k-1,\,k+1,\,\cdots,\,j$$ such that
$$
d_if_\ell =d_{\ell -1}f_i \text{\,\,\,if\,\,\,} i<\ell ,\,i\neq
k,\,\ell \neq k.$$ Then we have to find an $j$-simplex $$F\colon
K\times\Delta^m\rightarrow N$$ in $T_n(K)$ such that
$$d_iF=f_i.$$
By Hodgson \cite{Hodgson},

Let us think of $\Delta^j$ as a subpolyhedron of
$$
|\Lambda^k[j]|\times [0,1] \, .
$$
Then there is an embedded isotopy from the identity map of this
space to an embedding with image $\Delta^j$ (see e.g.,
\cite[p.\,295]{Hodgson}). Use this end of this isotopy to identify
$|\Lambda^k[j]|\times [0,1]$ with $\Delta^j$. Next, set
$$
N \,\, := \,\, (\bigcup_i {M_i})\times [0,1]\, .
$$
Then there we can define an extension
$$
F\colon K \times |\Lambda^k[j]|\times [0,1]\rightarrow N
$$
by the formula
$$F(x, (y, t))=(f_i(x, y), t)\, \text{\,\,\,if\,\,\,} (x, y)\in K\times
\Ds_i^{j{-}1}.$$
\end{proof}

\section{Stabilization}

\subsection*{The suspension map}
We consider the map
$$
E \: T_n(K) \to  T_{n+1}(K)
$$
defined by
$$
(M_\bullet,f_\bullet) \mapsto (M_\bullet \times [0,1],f'_\bullet)
$$
in which, for each $S$, the map  $f'_S$ is the composite
$$
\begin{CD}
K \times \Delta^S @> f_S >>  M_S = \frac{1}{2} \times M_S @>
\subset >> [0,1] \times M_S\, .
\end{CD}
$$
%(We use the right shift identification ${\mathbb R} \times
%{\mathbb R}^\infty = \mathbb R^\infty$ to consider the underlying
%manifold of $M_\bullet \times [0,1]$ as a submanifold of $\mathbb
% R^\infty$.)

Let
$$
T_\infty(K) \,\, := \,\, \lim_{n\to \infty} T_n(K)
$$
be the result of taking the colimit of the sequence defined by
$E$.

\begin{prop}\label{(stable is Kan)}
$T_{\infty}(K)$ is a Kan $\Ds$-set.
\end{prop}

\begin{proof} Follows from
that fact that $T_n(K)$ is Kan and the small object argument (i.e.,
the compactness of $\Delta[j]$).
\end{proof}

\subsection*{The tangent map}
A $j$-simplex $(N_\bullet,h_\bullet)$
of $T_n(K)$ has an underlying $(n{+}j)$-dimensional
manifold $N_{[j]}$
(which, by our conventions, is embedded in $\Bbb R^\infty$).
Let $BO_k$ be
the Grassmannian of $k$-planes in $\Bbb R^\infty$.
The usual Gauss map for the tangent bundle of $N_{[j]}$ gives a
$$
N_{[j]} \to BO_{n{+}j}
$$
We precompose this map with $h_{[j]}\: K \times \Delta^j \to N_{[j]}$
to get a map
$$
K\times \Delta^j \to BO_{n+j}\, .
$$
Let $BO$ be the direct limit of $BO_k$ with respect to the inclusions
$$
BO_k \subset BO_{k+1}
$$
defined by mapping a $k$-plane $V$ to the $k$-plane $\Bbb R \oplus V$.

Let
$$
 \text{maps}(K,BO)
$$
denote the total singular complex of the space of maps $K\to BO$,
considered as a $\Ds$-set. This satisfies the Kan condition.

Now let $j$ vary. Then the above construction assembles to a map
of $\Ds$-sets
$$
\tau \: T_n(K) \to \text{maps}(K,BO)\, .
$$
called the {\it tangent map}.

Note that with respect to $\tau$,
The stabilization map $E$ corresponds the self-equivalence of
$BO$ given by adding a copy of the real line.

\begin{thm} The tangent map $\tau$ defines a weak equivalence
of $\Ds$-sets
$$
T_\infty(K) \,\, \simeq \,\, \text{\rm map}(K, BO)\, .
$$
\end{thm}

\begin{proof} This is basically a result of
Wall \cite[Prop.\,\,5.1]{Wall1} adapted to block families. We will
sketch the argument in degree zero first, and then explain the
changes which are needed to adapt the argument to block families.

Observe  that if $L\overset \sim  \to K$ is a homotopy equivalence
of cell complexes, then  composition
gives an induced weak equivalence
$$
T_n(K) \overset\sim\to T_n(L)\, .
$$
This allows us to replace $K$ by a suitable compact framed manifold
within its homotopy type. Hence, without loss in generality,
it is enough to prove the result when $K$ itself is a compact
smooth manifold.

The next point is that a map $K \to BO$ is represented by a
{\it smooth} vector bundle $\xi\:E(\xi) \to K$, equipped with
a fiberwise inner product. The unit disk bundle $D(\xi)$ is then
a smooth manifold with boundary $D(\xi|\partial K) \cup S(\xi)$,
where $\xi|\partial K$ is the restriction of $\xi$ to $\partial K$
and $S(\xi)$ is the unit sphere bundle of $\xi$. The zero section
$z\: K \to D(\xi)$ then defines a thickening
$$
(D(\xi),z)\, .
$$
The tangent bundle of $D(\xi)$ pulled back along $z$ is canonically
identified with $\xi$. Hence, $\tau(D(\xi),z) \cong \xi$. This shows that
the map $\tau$ is surjective on $\pi_0$.

To obtain surjectivity in higher degrees, the
argument is modified as follows: a homotopy class
in $x\in \pi_j(\text{map}(K,BO))$  is represented by
a $j$-simplex $K \times \Delta^j \to BO$ such that
the restriction to $K \times \partial \Delta^j$
is constant in the second factor, i.e., it is a composite
$$
\begin{CD}
K \times \partial \Delta^j @> \text{project} >> K @> \alpha >> BO
\end{CD}
$$
Furthermore, by the degree zero case,
 we can assume that $\alpha = \tau(N,h)$ for some thickening $(N,h)$.

Applying the same argument as in degree zero,
but now working relatively, we can construct a $j$-simplex
of $T_\infty(K)$ which is constant along the boundary (i.e., it
has the form $(N \times \sigma,h\times \text{id})$
for each face $\sigma$ of $\Delta^j$), and whose tangential
data coincides with $x$.
This $j$-simplex yields the desired lift of $x$ to $\pi_j(T_\infty(K))$.
This completes the proof of surjectivity.

We now discuss injectivity in degree zero.
Fix two $n$-thickenings $(N_0,h_0)$ and $(N_1,h_1)$. Let
$\xi_i$ be the effect of applying $\tau$ to $(N_i,h_i)$,
and assume that there is a one simplex
$$
\xi\:K \times \Delta^1 \to BO
$$
which restricts to $\xi_i$ on $K \times i$.

Then there is a simple homotopy equivalence
$$
f\:N_0 \overset\sim \to N_1
$$
and an isomorphism of stable tangent bundles
$$
\xi_0 \,\, \cong \,\, f^*\xi_1\, .
$$

Since we are stabilizing, we can assume that $n \gg k$. Then by
transversality and Smale-Hirsch theory
(\cite{Hirsch},\cite{Smale}), $f\:N_0 \to N_1$ is homotopic to an
embedding into the interior of $N_1$, which also has the property
that $f^*(\tau_{N_1})$ is bundle isomorphic to $\tau_{N_0}$. Call
this embedding
$$
e\: N_0 \overset\subset\to N_1 \, .
$$
Then
$$
W = \text{closure}(N_1 - e(N_0))
$$
is an $s$-cobordism between $\partial N_1$ and $e(\partial N_0)$.
By the $s$-cobordism theorem, $W$ is diffeomorphic rel
 $e(\partial N_0)$ to a collar manifold $e(\partial N_0)\times [0,1]$

Then the manifold
$$
N_0 \cup W \cup N_1
$$
is diffeomorphic to the boundary of $N_{01}:= N_1 \times \Delta^1$
in which $N_0$ is identified with $N_1 \times 0$ and
$N_1$ is identified with $N_1 \times 1$. This gives a $3$-ad
$N_\bullet$ whose faces are $N_0$ and $N_1$.

Furthermore, there is
an evident simple block equivalence
$$
h_\bullet \:K \times \Delta^1 \to N_\bullet
$$
which extends the maps $h_i$. The pair $(N_\bullet,h_\bullet)$
gives a $1$-simplex of $T_\infty(K)$ whose faces
are $(N_i,h_i)$.  This shows injectivity in degree zero.

The case of injectivity in higher degrees is similar, and will
be left to the reader.
\end{proof}

\section{Compression}\label{compression}
%\resettheoremcounters

In the previous section we defined the suspension map
$$
E\: T_n(K)\rightarrow T_{n+1}(K).
$$
In this section we measure the deviation of $E$
from being surjective on path components.

\subsection*{Compression and sectioning}

\begin{defn} An $(n{+}1)$-thickening
$\alpha = (N,f)$  {\it compresses} if there is an $n$-thickening
$\beta$ and a $1$-simplex $\sigma$ of $T_{n+1}(K)$ such that
$$
d_0\sigma = \alpha \qquad \text{ and } \qquad d_1\sigma =
E(\beta)\, .
$$
\end{defn}

\begin{defn} The {\it fiberwise suspension} of a map of spaces
$E\to B$ is the map $S_B E \to B$ in which
$$
S_B E \,\, = \,\, B \times 0 \cup E\times [0,1]\cup B\times 1
$$
is the evident double mapping cylinder.
\end{defn}

When $E \to B$ is a (Hurewicz) fibration then so is $S_B E\to B$
(Str{\o}m \cite[p.\ 436]{Strom}) There are also two preferred
sections
$$
s_-,s_+ \: B \to S_B E\, .
$$

Although the following result is classical (cf. Larmore
\cite[Th.\,\,4.3]{Larmore}, see also  Becker\cite{Becker}), we
include a proof for the sake of completeness.

\begin{prop} \label{primitive_euler} Assume
that $E \to B$ is a fibration. If $E \to B$ admits a section then
$s_-$ and $s_+$ are section homotopic.

Conversely, let $E \to B$ be $(r{+}1)$-connected
 and assume that $B$ has the homotopy type of
a complex of dimension $\le 2r{+}1$. Furthermore, assume that
$s_-$ and $s_+$ are section homotopic. Then $E \to B$ admits a
section.
\end{prop}

\begin{proof}  A section $B \to E$ can be fiberwise
suspended over $B$ to give a map
$$
B \times [0,1] = S_B B \to S_B E\, .
$$
This map is the desired homotopy from $s_-$ to $s_+$.

Conversely, the commutative square
$$
\begin{CD}
E @>>> B \\
@VVV @VVs_+ V \\
B @>> s_- > S_BE
\end{CD}
$$
is homotopy cocartesian. By the Blakers-Massey theorem (see e.g.,
\cite[p.\,\,8]{Goodwillie}), it is also $(2r{+}1)$-cartesian. This
implies that the map from $E$ into the homotopy pullback of
$s_\pm$ (given by choosing a homotopy from $s_-$ to $s_+$) is
$(2r{+}1)$-connected. The conclusion follows.
\end{proof}

The next lemma, shows that homotopy and section homotopy coincide.
Although it will not be needed, we include it for the sake of
clarification.

\begin{lem}[Crabb-James {\cite[p.\,\,29]{Crabb-James}}] Let $E \to B$ be a Hurewicz
fibration and suppose that $s,t\: B \to E$ are sections. Then $s$
and $t$ are homotopic if and only if they are section homotopic.
\end{lem}

%\begin{proof} The `only if' part is clear. To prove the `if' part,
% let $\text{map}(B,E)$ be the space of maps from $B$ to $E$. Let
% $p\:E\to B$ denote the fibration. Then $p$ induces a map
% $$
% p_*\: \text{map}(B,E) \to \text{map}(B,B)
% $$
% which is again a fibration, and the fiber over the identity map of
% is clearly the space of sections of $p$. Composition with $s$
% induces a section of $p_*$. In particular, if $\text{sec}(p)$
% denotes the section space of $p$, then the fiber inclusion
% $$
% \text{sec}(p) \subset \text{map}(B,E)
% $$
% induces an injection on homotopy---in particular on path components.
% This completes the proof.
%\end{proof}

Returning to the original problem of compressing $\alpha$, we let
$E_\alpha$ be the homotopy pullback of the diagram
$$
\begin{CD}
K @> f >\sim > N @< \supset << \partial N
\end{CD}
$$
Then $E_\alpha \to K$ is a fibration.

\begin{thm} \label{prototype} If $\alpha$ compresses, then
the preferred sections
$$
s_\pm\: K \to S_K E_\alpha
$$
are homotopic.

Conversely, Assume $3k+2 \le 2n$, where $K$ has the homotopy type
of a complex of dimension $\le k$. Then the existence of a
homotopy between $s_-$ and $s_+$ implies that $\alpha$ compresses.
\end{thm}

\begin{proof} To prove the first part, first assume that
$\alpha = E(\beta)$, where $\beta = (M,g)$. Then $N = M \times
[0,1]$ and $M \times 0 \subset \partial N$. The map $g\: K \to
M\times 0$ followed by that inclusion into $\partial N$ then gives
rise to the section of $E_\alpha \to K$ in this instance. By the first part
of \ref{primitive_euler}, we get a homotopy between $s_-$ and $s_+$ in this
instance.

The more general case when $\alpha$ compresses uses the
$s$-cobordism between $\partial N$ and $\partial (M \times [0,1])$
which is provided  by the $1$-simplex $\sigma$.

We now prove the converse. Our range assumptions imply that the
map $E_\alpha \to K$ is $(n{-}k)$-connected. Hence, by
the second part \ref{primitive_euler},
 the existence of a homotopy between $s_+$
and $s_-$ enables us to conclude that the simple homotopy
equivalence $f\: K \to N$ factors through $\partial N$ up to
homotopy. Let $f_1\: K \to
\partial N$ be a choice of such a factorization.

The rest of the proof appeals to Wall's original argument which we
will now sketch. The map $K \to
\partial N$ is $(n-k-1)$-connected (since $\partial N \to N$ is
$(n-k)$-connected). Our assumptions guarantee that
$$
n - k - 1 \ge 2k - n + 1
$$
so by the Stallings-Wall embedding theorem (see
\ref{stallings-wall} below) one can find an {\it embedding up to
homotopy} of the map, i.e., there is a codimension zero compact
submanifold $M \subset
\partial N$
(satisfying the $\pi$-$\pi$-condition)
and a simple homotopy equivalence $g\:K \overset\sim\to M$ such
that the composite
$$
K \overset g \to M \subset \partial N
$$
is homotopic to the given map. Then $(M,g)$ is an $n$-thickening
such that $(M,g)$ and $\alpha$ are the faces of a one simplex of
$T_{n+1}(K)$ where the underlying manifold of the one simplex is
$$
N \cup_{M\times 0} M \times [0,1] \, .
$$
\end{proof}

\begin{thm}[Stallings {\cite[p.\ 5]{Stallings}}, Wall {\cite[p.\ 76]{Wall1}}]\label{stallings-wall}
Let $K$ be a connected finite complex. Let $M$ be a manifold of
dimension $n$ with boundary. Assume $\dim K\le k \le n {-}3$. If
$f\colon K\to M$ is $(2k {-} n {+} 1)$-connected, then $f$ embeds
up to homotopy.
\end{thm}

\subsection*{The Euler class}

%If $\alpha$ compresses,
%it follows that $f^*(\tau_N)$ is isomorphic
%to $\xi \oplus \epsilon$, for
%some $\xi$ a rank $n$ vector bundle over $K$.

%\begin{ass} There is fixed chosen isomorphism of vector bundles
%$$
%\xi \oplus \epsilon \,\, \cong \,\, f^*(\tau_N)
%$$
%\end{ass}

%With respect to this assumption, we will define a secondary obstruction
%to deforming $\alpha$ to an element of the form $E(\beta)$.

%For $E \to B$ a fibration, let
%$$
%\sec(E\to B)
%$$
%denote its space of sections.

For an $(n{+}1)$-thickening $\alpha = (N,f)$, set
$$
\xi \,\, := \,\, f^*(\tau_N)\, .
$$
For $x\in K$ let
$$S^\xi_x$$ denote the one point compactification
of the fiber of $\xi$ at $x$. For $(x,y) \in K \times K$, let
$$
\Omega_x^y K
$$
denote the space of paths in $K$ which begin at $x$ and end at
$y$.

There is a fibration
$$
{\cal E}_\xi\to K \times K
$$
whose fiber at $(x,y) \in K$ is the space
$$
S^\xi_x \smsh (\Omega_x^y K)_+\, .
$$
Explicitly, if we equip $\xi$ with a fiberwise inner product, then
${\cal E}_\xi$ can be described by the space of pairs
$(v,\lambda)$ in which $\lambda$ is a path in $K$ and $v$ is a
point in the unit disk of $\xi$ at $\lambda(0)$, subject to the
identification $(v,\lambda) \sim (w,\beta)$ for all $v$ and $w$ of
length one. With this description the map ${\cal E}_\xi\to K
\times K$ is given by $(v,\lambda) \mapsto
(\lambda(0),\lambda(1))$. This fibration comes equipped with a
preferred section $s\:K \times K \to {\cal E}_\xi$.

\begin{defn} Let $E \to B$ be a (Hurewicz) fibration
equipped with a preferred section $s\:B \to E$. Let
$$
Q_\bullet E \to B
$$
be the effect of applying the stable homotopy functor $Q =
\Omega^\infty\Sigma^\infty$ fiberwise to $E \to B$ (\cite{Crabb-James}).

The {\it stable section space} of $E \to B$ is the space of
sections of $Q_\bullet E \to B$. It will be denoted by
$$
\sec^{\rm st}(E \to B)\, .
$$
\end{defn}

\begin{rem} The stable section space has the structure of
an  infinite loop
space. In particular, its set of path components is an abelian group.
\end{rem}

Recall that $E_\alpha \to K$ is the fibration given by taking the
homotopy pullback of $\partial N \to N$ along $f\: K \to N$. The
fiberwise suspension $S_KE_\alpha \to K$ has two preferred
sections $s_\pm$. Our convention will be to use  $s_-$ to form its
stable section space.

\begin{thm}\label{section_ident}
There is a weak equivalence of stable section spaces
$$
\sec^{\rm st}(S_K E_\alpha \to K) \,\, \simeq \,\, \sec^{\rm
st}({\cal E}_\xi \to K \times K)
$$
which can be made canonical up to  contractible choice.
\end{thm}

The proof of the theorem is deferred until the next section.

\begin{rem} The  right hand side only
depends on {\it stable isomorphism type} of
the bundle $\xi$, not on the choice of thickening
$\alpha = (N,f)$.
\end{rem}

\begin{defn} The {\it Euler class}
$$
e(\alpha) \in \pi_0(\sec^{\rm st}({\cal E}_\xi \to K \times K))
$$
is the element which corresponds to the section  $s_+\in \sec^{\rm
st}(S_K E_\alpha \to K)$ via the weak equivalence of
\ref{section_ident}.
\end{defn}

\begin{lem}\label{same}
Suppose $p\:E\rightarrow B$ is a fibration equipped with section
$B \to E$. Assume that $p$ is $(r{+}1)$-connected and that
$B$ has the homotopy type of a cell complex of dimension
$\leq b$. Then the evident map
$$\sec(E\rightarrow B)\rightarrow
\sec^\st(E\rightarrow B)$$ is $(2r{+}1{-}b)$-connected.
\end{lem}

\begin{proof} Let $F$ denote any fiber of $p$. Then
$F$ is an $r$-connected based space and the map
$$
F \to QF
$$
is $(2r{+}1)$-connected by the Freudenthal suspension theorem.
The long exact homotopy sequence and the 5-lemma then imply that the map
$$
E \to Q_\bullet E
$$
is also $(2r{+}1)$-connected. The induced map of section spaces is
then $(2r{+}1{-}b)$-connected by obstruction theory.
\end{proof}

\begin{cor} If $\alpha$ compresses, then $e(\alpha)$ is trivial.
Conversely, if $3k{+}2 \le 2n$ and $e(\alpha)$ is trivial, then
$\alpha$ compresses.
\end{cor}

We note that if $n \geq 2k$, then the section space
in lemma \ref{same} is path connected (by obstruction theory). Then
$H(\alpha)=0$ for any vertex $\alpha$ in $T_{n{+}1}(K)$, and it
follows that $E$ is surjective in this case.

\subsection*{Equivariant homotopy theory}

Let $B$ be a connected space. Define
$$\Top_{/B}$$ to be the category of spaces ``over $B$.''
Specifically, an {\it object}
is a space $X$ and a choice of (structure) map $X\to B$. A {\it morphism} is a
map of spaces which is compatible with the structure maps. Call a
morphism a {\it weak equivalence}, if it is a weak homotopy
equivalence of underlying  spaces.
%{\it fibration} or {\it cofibration} if it respectively is so
% when considered in $\Top$ by
%means of the forgetful functor $\Top_{/B}\to \Top$. $\Top_{/B}$
%endowed with these three families of maps is a {\it model
%category} [Klein].
It can be shown that this notion of weak equivalence arises
from a Quillen model structure on $\Top_{/B}$, but the full strength
of this
fact will not be needed.
\medskip

Next, let $G$ be a topological group which is cofibrant
a considered as a topological space.
Define $$\C_G$$ to be the category of (left) $G$-spaces and $G$-equivariant
maps. A {\it weak equivalence} in this instance
is deemed to be a morphism
whose underlying map of (unequivariant) spaces is a weak homotopy
equivalence. Again, this arises from a Quillen model structure.
\medskip

In each of these cases, let the homotopy category be formed in the
usual way by formally inverting the weak equivalences.
\medskip

\subsection*{The correspondence}
The prototype idea is that any connected based space is naturally
the weak homotopy type of a classifying space of a topological
group, the latter being a suitable model for the loop space.

Specifically,
if $B$ is a connected based space, then a construction of Kan
gives a natural weak equivalence
$$
B \,\, \simeq \,\, BG
$$
where $G$ is a topological group object in the category of
compactly generated spaces.
The definition of $G$ is as follows:
Let $G.(B)$ denote the Kan loop group of the total singular
complex of $B$. Then $G$ is taken to be the geometric realization
of $G.(B)$\, .

Define a functor
$$
\C_G \to \Top_{/BG}
$$
by $X \to X \times_G EG$ (the Borel construction).

\begin{prop}[cf. {\cite[Th.\ 2.1 and Cor.\ 2.5]{Dror-Dwyer-Kan}}] This functor
can be derived to give an equivalence of homotopy categories
$$
\text{\rm ho\,}\Top_{/BG}\,\, \simeq \,\, \text{\rm ho\,}\C_G\, .
$$

\end{prop}

It will be helpful to have
functor which gives rise to the inverse
equivalence.

To this end, let $X \to BG$ be a map of spaces. Define its {\it thick
homotopy fiber} by
$$
F := \text{pullback}(EG \to BG \leftarrow X)
$$
in which $EG \to BG$ is the universal principal $G$-bundle. Then
$F \subset EG \times X$ is a subspace, $G$ acts (freely and on the
left) on the first factor $EG$ and therefore on $EG \times X$
using the trivial action on $X$. It is readily verified that this
restricts to an action of $G$ on $F$. The assignment $X \mapsto F$
defines a functor $\Top_{/BG}\to \C_G$ which induces the inverse
equivalence on homotopy categories.
\medskip

Here is a little more detail which contains the crux of
the statement:
as $EG$ is a contractible space, we see that $F$ is a
model for the (usual) homotopy fiber of $X\to BG$ at the preferred
basepoint. The following shows how to recover $X\to BG$ up to
homotopy from its thick homotopy fiber:

\begin{lem} Assume that $X$ is a cofibrant space.
Then in the homotopy category of spaces over $BG$ there is a
canonical weak equivalence between the Borel construction
$$
EG \times_G F \to BG
$$
and $X \to BG$.
\end{lem}

\begin{proof}
The displayed map factors as
$$
EG \times_G F \to X \to BG
$$
In which the first map is obtained from the equivariant weak
equivalence
$$
EG \times F \to \text{pt} \times F
$$
by taking $G$-orbits (we are using the observation that $X$ is
obtained from $F$ by taking $G$ orbits). Since $F$  has the
structure of a free $G$-CW complex (by covering homotopy
property), both source and target of the above equivariant weak
equivalence
 are fibrant and cofibrant. This implies that
the map is an equivariant homotopy equivalence. Hence the map of
orbit spaces is a homotopy equivalence.
\end{proof}

\subsection*{Dictionary} The following table
gives the correspondence between various constructions
in the categories $\Top_{/B}$ and $G$-$\Top$.

$$
\vbox{ \offinterlineskip \tabskip=2pt \halign{ \strut # & \vrule #
& \hfil # \hfil & \vrule # & \hfil # \hfil & \vrule # & \hfil #
\hfil & \vrule #\cr \omit& \multispan{7}{\hrulefill}\cr & & {\bf
Fibration} $E \to B= BG$ && & & $F=$ {\bf thick homotopy fiber} &
\cr \omit& \multispan{7}{\hrulefill}\cr & & section $B \to E$   &
& & & $F=$ based $G$-space &\cr & & & & & & &\cr \omit&
\multispan{7}{\hrulefill}\cr & & $S_B E \to B$ & & & & $SF=$
unreduced suspension &\cr & & & & & & &\cr \omit&
\multispan{7}{\hrulefill}\cr & & $\text{sec}(E \to B)$ & & & &
$F^{hG}$ = homotopy fixed points &\cr & & & & & & &\cr \omit&
\multispan{7}{\hrulefill}\cr & & $\text{sec}^{\rm st}(S_B E \to
B)$ & & & & $(QSF)^{hG}$ &\cr & & & & & & &\cr \omit&
\multispan{7}{\hrulefill}\cr }}
$$

\subsection*{Naive $G$-spectra}
A {\it (naive) $G$-spectrum} $E$
consists of based (left) $G$-spaces
$E_i$ for $i \ge 0$, and equivariant based maps $\Sigma E_i \to
E_{i{+}1}$ (by convention, $G$ acts  trivially on the suspension
coordinate of $\Sigma E_i$). A {\it morphism} $E \to E'$ of
$G$-spectra consists of equivariant
 based maps $E_i \to E'_i$ which are compatible with
the structure maps.  A {\it weak equivalence}  of $G$-spectra is a
morphism that yields an isomorphism on homotopy groups. A
$G$-spectrum $E$ is said to be  {\it fibrant} if the adjoint maps
$E_i \to \Omega E_{i{+}1}$ are weak homotopy equivalences.
Up to a
natural weak equivalence,
one can always approximate a $G$-spectrum by one which is fibrant
(this procedure is called {\it fibrant replacement}).
The {\it zeroth space} $E_0$ of a fibrant spectrum
$E$ is often denoted $\Omega^\infty E$. If $E$ isn't fibrant,
we take $\Omega^\infty E$ to the zeroth space its fibrant replacement.

If $X$ is a based $G$-space, then its {\it suspension spectrum}
$\Sigma^{\infty} X$ is a $G$-spectrum with $j$-th space $Q(S^j
\smsh X)$, where $Q = \Omega^{\infty}\Sigma^{\infty}$ is the
stable homotopy functor (here $G$ acts trivially on the suspension
coordinates). We use the notation
$$
S[G]
$$
for the suspension spectrum of $G_+$ considered as a
$(G{\times}G)$-spectrum (the action on $G_+$ is given by left
multiplication with respect to the first factor of $G{\times} G$
and  right multiplication composed with the involution $g \mapsto
g^{-1}$ on the second factor.

%If $E$ is a  $G$-spectrum then the {\it homotopy orbit spectrum}
%$E_{hG}$ is
%$$
% E \smsh_G EG_+\, ,
%$$
%where $EG$ the free contractible $G$-space (arising from the bar
%construction), and $EG_+$ is the effect of adding a basepoint to
%$EG$.

Let $E$ be a  fibrant $G$-spectrum. Then its
{\it homotopy fixed point spectrum} $E^{hG}$ is
$$
\text{\rm map}_G(EG_+,E)
$$
where the $j$-th space of the latter is given by the mapping space
of equivariant maps $EG_+ \to E_j$. Here $EG$ the free
contractible $G$-space (arising from the bar construction), and
$EG_+$ is the effect of adding a basepoint to $EG$.

\subsection*{The dualizing spectrum}
Using the left factor action $G$ on $S[G]$ we can  form the
homotopy fixed point spectrum
$$
D_G := S[G]^{hG}\,
$$
(cf.\ Klein \cite{Klein}).
This is called the {\it dualizing spectrum of $G$.}
The right factor action of $G$ on
$S[G]$ gives $D_G$ the structure of a $G$-spectrum.

\section{Proof of Theorem $\ref{section_ident}$}

\subsection*{The fibration ${\cal E}_\xi \to K {\times} K$}
Recall that the fiber of ${\cal E}_\xi \to K {\times} K$ at
$(x,y)$ is identified with
$$
S^\xi_x \smsh (\Omega_x^y K)_+
$$
Let $* \in K$ denote the base point. In particular,
 the fiber at $(*,*)$ is given
by
$$
S^\xi_* \smsh (\Omega K)_+ \, ,
$$
in which $\Omega K$ is the based loop space of $K$. We will give
an equivariant model for this space.

In what follows by slight abuse of notation, we identify $K$ by
$BG$. Then we think of $\xi$ is a vector bundle over $BG$; let
$$
S^\xi := {\widetilde D(\xi)}/{\widetilde S(\xi)}\, ,
$$
${\widetilde D(\xi)}$ (resp.\ ${\widetilde S(\xi)}$)
denotes the thick homotopy fiber of
$D(\xi) \to BG$ (resp.\ $S(\xi) \to BG$).
Then $S^\xi$ is a $G$-equivariant
model for $S^\xi_*$. Extend this to an  action of $G \times G$ by
letting the second factor act trivially. Let $G^{\text{ad}}$
denote $G$ with the action of $G \times G$ given by
$$
(g,h)\cdot x := gxh^{-1}\, .
$$
Let $G^{\rm ad}_+$ be the effect of adding a disjoint basepoint to
$G^{\rm ad}$.

Lastly, set
$$
S^\xi[G] \,\, := S^\xi \smsh G^{\rm ad}_+ \, ,
$$
and give this the diagonal $G\times G$ action.

With respect to these conventions we obtain

\begin{lem}\label{reduction_1}
In the category of spaces over $BG \times BG$, there is a weak
equivalence
$$
{\cal E}_\xi \,\, \simeq \,\, (EG\times EG)\times_{G\times G}
S^\xi[G]
$$
where $G \times G$ acts diagonally on $S^\xi \smsh G^{\text{\rm
ad}}_+$. In particular, there is a weak equivalence of infinite
loop spaces
$$
\sec^{\rm st}({\cal E}_\xi \to K \times K) \,\, \simeq \,\, (Q
S^\xi[G])^{hG\times G}\, .
$$
\end{lem}

\begin{proof} The fibration
$$
 {\cal E}_\xi \to K \times K
$$
is the {\it fiberwise smash product} of two fibrations over $K
\times K$. The first of these fibrations is given as follows:
consider the vector bundle $\xi$ over $K$, apply fiberwise one
point compactification and take the resulting pullback along first
factor projection $K \times K \to K$. This results in the
fibration $E_1 \to K \times K$ whose fiber at $(x,y)$ is
$S^\xi_x$.

The second fibration over $K \times K$ can be described as
follows: let $K^{[0,1]}$ be the space of paths in $K$, and let
$$
K^{[0,1]} \to K \times K
$$
be the fibration which evaluates a path at its endpoints. The
fiber of this fibration at $(x,y)$ is then precisely $\Omega_x^y
K$. Consequently, if we add a basepoint to each fiber, we obtain a
fibration $E_2 \to K \times K$ whose fiber is $(\Omega_x^y K)_+$.

Clearly, we have a homeomorphism over $K \times K$
$$
 {\cal E}_\xi  \,\, \cong \,\, E_1 \smsh_{K\times K} E_2\, ,
$$
where the right side denotes the fiberwise smash product of $E_1$
and $E_2$.

The map $K^{[0,1]} \to K \times K$ is the effect of converting the
diagonal $K \to K \times K$ into a fibration. In terms of the
identification $K \simeq BG$, this diagonal map can be identified
with the Borel construction of $G\times G$ acting on
$G^{\text{ad}}$. The reason for this is that one can consider
$$
(EG \times EG) \times_{G {\times} G} G^{\text{\rm ad}}
$$
as being formed in 2 steps:
\begin{enumerate}
\item  take the Borel construction on $G^{\text{\rm ad}}$ with
respect to the subgroup $G \times 1 \subset G\times G$, and then
\item  the Borel construction with respect to $1 \times G$.
\end{enumerate}
The result of the first step gives the $G$-space $EG$ equipped
with the usual map to $BG$. Consequently, application of the
second step gives the evident map
$$
EG \times_G EG \to BG \times BG \, .
$$
The total space of this fibration is identified with $BG$ by
taking the $G$-orbits of the the diagonal map $EG \to EG \times
EG$. With respect to this identification, the displayed map
coincides with the diagonal of $BG$.

The proof is completed by observing that the fiberwise smash
product, corresponds via the dictionary to the smash product of
equivariant spaces.
\end{proof}

\subsection*{The exponential law} Let $E$ be a (fibrant)
$(G\times H)$-spectrum. Then we have an equivalence of homotopy
fixed point spectrum in two steps
$$
E^{hG\times H} \,\, \simeq  \,\, (E^{hH})^{hG} \, .
$$
If we apply this reasoning to $G\times G$ acting on the spectrum $
(\Sigma^\infty S^\xi[G]) $, we obtain, after a some trivial
rewriting

\begin{lem} \label{reduction_2}
There is a weak equivalence of spectra
$$
(\Sigma^\infty S^\xi[G])^{hG\times G} \simeq (S^\xi \smsh
D_G)^{hG}
$$
\end{lem}

Combining this with \ref{reduction_1} then yields

\begin{cor} \label{combined} There is a weak equivalence of infinite loop
spaces
$$
\sec^{\rm st}({\cal E}_\xi \to K \times K) \,\, \simeq \,\,
\Omega^\infty (S^\xi \smsh D_G)^{hG}
$$
\end{cor}

\subsection*{Identification of $E_\alpha \to K$}
Consider now the fibration
$$
E_\alpha \to K\, .
$$
Make the identification $K \simeq BG$ and let $F$ denote the thick
homotopy fiber. (Recall that $E_\alpha \to K$ is identified with
the inclusion $\partial N \to N$ up to homotopy.) The following is
essentially \cite[Prop.\ 10.5]{Klein}.

\begin{prop}[Klein]
\label{reduction_0} There is an equivariant weak equivalence
$$
S^\xi \smsh D_G \, \, \simeq \,\, \Sigma^\infty SF\, .
$$
\end{prop}

\begin{rem}
The case in which $\xi$ is the trivial fibration is one of the
main theorems of \cite{Klein} which relates the Spivak fiber (the
right side) to equivariant homotopy theory (the left side).

The above weak equivalence is defined by means of a certain kind
of equivariant duality map. In particular, the equivalence is
given up to contractible choice.
\end{rem}

%We are now in position to prove \ref{section_ident}.

\begin{proof}[Proof of Theorem \ref{section_ident}]

\begin{align*}
\sec^{\rm st}(S_K(E_\alpha \to K)) \quad & \simeq \quad
(QSF)^{hG}  \qquad \text{ by the dictionary,} \\
& \simeq \quad \Omega^\infty (S^\xi \smsh D_G)^{hG}
\qquad \text{ by } \ref{reduction_0}, \\
& \simeq \quad \sec^{\st}({\cal E}_\xi \to K \times K) \quad
\text{ by } \ref{combined}.
\end{align*}
\end{proof}

\section{The homotopy fibers of $E$}
%\resettheoremcounters

In this section we identify the homotopy
fibers of the suspension map in a stable range.

\subsection*{Abstract description}

Let $\alpha = (N,f)$ be a 0-simplex of $T_{n{+}1}(K)$.
We give the description of the homotopy fiber of the suspension
map $E$ at the basepoint $\alpha$ (cf.\ Chapter 1):

\begin{defn} Let $F_n(\alpha)$ be the $\Ds$-set in which a $j$-simplex
is specified by a pair
$$
(\beta,\sigma)
$$
in which $\beta\in T_n(K)$ is a $j$-simplex and $\sigma\:\Delta[j]
\times \Delta[1] \to T_{n+1}(K)$ is a map of $\Ds$-sets such that
$\sigma$ restricted to $\Delta[j]\times 0$ represents
$\alpha\times id_{\Delta^{j}}$ and $\sigma$ restricted to
$\Delta[j] \times 1$ represents $E(\beta)$.
\end{defn}

\subsection*{Formulation}

The last chapter gave criteria for deciding when $F_n(\alpha)$ is
non-empty in the metastable range. The goal of this chapter is to
identify $F_n(\alpha)$ in the metastable range. We will define a map
$$
{\cal H}\: |F_n(\alpha)| \to
\Omega \, \sec^{\st}({\cal E}_\xi \to K \times K)\, ,
$$
where the space on the right hand side is the loop space of the
stable  section space which appeared in the previous chapter. We will
then show that this map is highly connected:

\begin{thm} \label{euler-map-thm} Assume that $K$ has
the homotopy type of a cell complex of dimension $\le k$ and that
$\alpha$ admits a compression.

Then there is $(2n{-}3k{-}2)$-connected map
$$
{\cal H}\: |F_n(\alpha)| \to
\Omega \, \sec^{\st}({\cal E}_\xi \to K \times K)\, .
$$
\end{thm}
This will be the main result of the section.

\subsection*{Section spaces revisited}

Let $p\:E\to B$ be a (Hurewicz) fibration.
Consider as before the fiberwise suspended fibration
$$
S_B E \to B\, .
$$
Denote its fiber at $b\in B$ by $SF_b$. This space
comes equipped with two points given by the
south and north poles of the suspension.

\begin{defn} Let
$$
\Omega_B^\pm S_B E \to B
$$
be the fibration in which the fiber at $b\in B$ consists
of a path from the south to the north pole of $SF_b$.
\end{defn}

Explicitly, $\Omega_B^\pm S_B E$ is the homotopy pullback of the
diagram
$$
\begin{CD}
B @> s_- >> S_B E@< s_+ << B\, .
\end{CD}
$$
There is an evident map
$$
E \to \Omega_B^\pm S_B E
$$
which is $(2r{+}1)$-connected whenever $E\to B$ is
$(r{+}1)$-connected (this follows from the Blakers-Massey theorem).
It follows that the map of section spaces
$$
\text{\rm sec}(E\to B) \to
\text{\rm sec}(\Omega_B^\pm S_B E\to B)
$$
is $(2r{+}1{-}b)$-connected whenever $B$ has the homotopy type of
a cell complex of dimension $\le b$.

Now if the section space on the right happens to come equipped
with a basepoint, then any other section may be combined with the
basepoint section to form a {\it based loop of sections} (the
reason is that two paths between two points can be regarded as a
based loop at one of the points.) Summarizing

\begin{lem}\label{lem3.1} Given a choice of
basepoint in $\text{\rm sec}(\Omega_B^\pm S_B E\to B)$,
there is a weak equivalence of spaces
$$
\text{\rm sec}(\Omega_B^\pm S_B E\to B)\,\,
\simeq \,\,
\Omega \text{\rm sec}(S_B E\to B)\, .
$$
\end{lem}

\begin{cor}\label{cor3.1} Once a basepoint in
$\text{\rm sec}(\Omega_B^\pm S_B E\to B)$ is chosen, there is a
$(2r{+}1-b)$-connected map
$$
\text{\rm sec}(E\to B) \to \Omega \text{\rm sec}(S_B E\to B)\, .
$$
\end{cor}

We now give a convenient model for section spaces.

\begin{defn} Let $p\: E\to B$ be a fibration.
Then
%$$
%\text{\rm sec}{\bold .}(p)$\text{\rm sec}{\bold .}(E\to B)$
%$$
$$
\text{\rm sec}{\bold .}(E\to B)
$$
is defined to be the $\Ds$-set in which a $j$-simplex consists of
$$
(p',s)
$$
in which
\begin{itemize}
\item
$$
p'\:E' \to B \times \Delta^j \times [0,1]
$$
is a fibration,
\item the restriction of $p'$ to
$B\times \Delta^j \times 0$ coincides with
$$
p \times {\rm id}\:E \times \Delta^j \to B \times \Delta^j
$$
(where $p\: E\to B$ is the given fibration), and
\item $s$ is a section of $p'$ along the subspace
$B \times \Delta^j \times 1$.
\end{itemize}
The face operators are induced in the evident way from the faces
of $\Delta^j$.
\end{defn}

\begin{lem}\label{lem3.2} Assume that $B$ is a cell complex.
Then $\text{\rm sec}{\bold .}(E\to B)$
satisfies the Kan condition.
\end{lem}

\begin{proof}(Sketch).
To keep notation in the proof consistent, we rename
$\sec{\bold .}(E\to B)$ by
$$
\sec{\bold .}(p) \, .
$$
Let ${\cal F}(E\to B)$ be the $\Ds$-set
whose $j$-simplices are fibrations over $B\times \Delta^j$.
It is straightforward to check that ${\cal F}(E\to B)$
satisfies the Kan condition.

Let $P_p(E\to B)$ be the $\Ds$-set of paths
in  ${\cal F}(E\to B)$ with initial point $p\:E\to B$. Then
$P_p(E\to B)\to {\cal F}(E\to B)$ is a Kan fibration. Hence
$P_p(E\to B)$ is also a Kan $\Ds$-set.

The projection $(p',s) \to p'$ defines a
map
$$
\sec.(p) \to P_p(E\to B)\, .
$$
As $P_p(E\to B)$ is Kan, it is sufficient to verify that
this map satisfies the Kan condition. But this follows
directly from the homotopy extension principle.
\end{proof}

\begin{lem} \label{section_equiv}
Assume that $B$ is a cell complex.
Then there is a preferred weak equivalence of spaces
$$
|\text{\rm sec}{\bold .}(E\to B)| \,\, \simeq \,\, \text{\rm
sec}(E\to B)
$$
where the right side denotes the space of sections of $p\:E\to B$.
\end{lem}

\begin{proof} Let the notation be as in the previous proof.
Then as above we have a Kan fibration
$$
\sec{\bold .}(p)\to P_p(E\to  B) \, .
$$
 The fiber over the basepoint
is just
$$
S{\bold .}(\text{sections}(p))
$$
= the total singular complex of the space of
sections of $p$. This has the homotopy type of
the space of sections of $p$.

As $P_p(E\to B)$ is contractible, it follows that the inclusion
from fiber to total space is a weak equivalence. The conclusion follows.
\end{proof}

Combining  the above
with \ref{same}, we obtain

\begin{cor}\label{cor3.2} Assume  $E \to B$ is an $(r{+}1)$-connected
 fibration and $B$ is a cell complex
of dimension $\le b$. Assume further that the fibration
$\Omega^\pm_B S_B E\to B$ comes equipped with a section.

Then there is a $(2r{-}b)$-connected (weak) map
$$
|\sec{\bold .}(E\to B)| \to \Omega \sec^{\rm st}(S_BE\to B)\, .
$$
\end{cor}

\subsection*{Defining the map ${\cal H}$}

The first step is to define a map of $\Ds$-sets
$$
\Gamma\:F_n(\alpha) \to \sec{\bold .}(E_\alpha \to K)\, ,
$$
where we recall that $E_\alpha \to K$ is the fibration
defined by taking the homotopy pullback
of
$$
K \overset f \to N \leftarrow \partial N\, .
$$
We first describe the map on 0-simplices. Recall that a zero
simplex of $F_n(\alpha)$ consists of a $1$-simplex
$$
\sigma = (N_\bullet,f_\bullet)
$$
of $T_{n+1}(K)$ such that $d_1 \sigma$ is of the form $E(\beta)$
for some $\beta = (M,g)$ a 0-simplex of $T_n(K)$ and $d_0 \sigma =
\alpha$. Recall that $E(\beta)$ is the thickening given by
$$
K \overset g \to  M = M \times \frac{1}{2} \subset M \times [0,1]\, .
$$
Then we have a fibration
$$
E_\sigma \to K \times \Delta^1
$$
given by taking the homotopy pullback of
$$
\begin{CD}
K \times \Delta^1 @>>> N_{01} @<\supset<< \partial_1 N_{01}  \, .
\end{CD}
$$
The fibration has the following properties:
\begin{itemize}
\item its restriction
to $K \times 0$ is identified with $E_\alpha \to K$, and
\item its restriction to $K \times 1$ comes equipped with a
preferred section.
\end{itemize}
(The second property uses the preferred homotopy from
$M \times \frac{1}{2} \subset M \times [0,1]$ into $M \times 0 \subset
\partial (M \times [0,1])$.)

Consequently, these data describe a zero simplex of $\sec{\bold
.}(E_\alpha \to K)$. The map on higher dimensional
simplices is given by essentially the same description (we omit the details).

In summary, what we have constructed is a map of
$\Ds$-sets
$$
\Gamma\:F_n(\alpha) \to \sec{\bold .}(E_\alpha \to K)\, .
$$

Let us now assume that $\alpha$ admits a compression. Then the fibration
$E_\alpha \to K$ admits a non-trivial section and so does
$\Omega^\pm_K S_K E_\alpha \to K$.

Then by corollary \ref{cor3.2}, there is a weak map of spaces
$$
|\sec{\bold .}(E_\alpha \to K)| \to \Omega \sec^{\rm st}(S_KE_\alpha \to K)
$$
By theorem \ref{section_ident}, we have a weak equivalence
$$
\Omega \sec^{\rm st}(S_KE_\alpha \to K) \simeq \Omega \sec^{\rm st}({\cal E}_\xi \to K \times K)\, .
$$
Assembling the above, and applying corollary \ref{cor3.2}, we get

\begin{cor} \label{blah} There is a (weak) map of spaces
$$
|\sec{\bold .}(E_\alpha \to K)| \to
\Omega \sec^{\rm st}({\cal E}_\xi \to K \times K)\, .
$$
If $K$ is a cell complex of dimension $\le k$
and $\alpha$ is an $(n{+}1)$-thickening, then this map
is $(2n-3k-2)$-connected.
\end{cor}

\begin{defn} Fix a compression of $\alpha$.
Then the map
$$
{\cal H}\: |F_n(\alpha)| \to \Omega \sec^{\rm st}({\cal E}_\xi \to K \times K)
$$
is given by precomposing the map of \ref{blah} with
the realization of the map
$$
\Gamma\:F_n(\alpha) \to \sec{\bold .}(E_\alpha \to K)
$$
\end{defn}

\subsection*{Connectivity of ${\cal H}$}

We are now ready to prove the main result of this chapter. By
\ref{blah} and the given assumptions, it will be sufficient to
prove

\begin{prop}\label{Gamma_connect}
The map
$$
\Gamma\: F_n(\alpha) \to \sec{\bold .}(E_\alpha \to K)
$$
is $(2n-3k-2)$-connected (after realization).
\end{prop}

\subsubsection*{Digression: Hodgson's embedding theorem}
Let $(K,L)$ be a cofibration pair such that $L$
is a finite complex
and $K$ is obtained from $L$ by attaching cells
of dimension $\le k$. In this instance, we write
$$
\dim(K,L) \le k \, .
$$
Let $N$ be a compact $n$-manifold.

\begin{defn} Fix a map of pairs
$$
f:= (f_K,f_L)\: (K,L) \to (N,\partial N)\, .
$$
Then a {\it (relative) embedding up to homotopy}
of $f$ in $(N,\partial N)$ consists of a pair
$$
(W,h)
$$
such that
\begin{itemize}
\item
$W\subset N$ is a compact $n$-manifold
equipped with
boundary decomposition
$$
\partial W = \partial_0 W \cup_{\partial_{01} W} \partial_1 W
$$
such that $W$ intersects $\partial N$ transversely and $W \cap
\partial N = \partial_0 W$. \item
$$
h := (h_K,h_L)\: (K,L) \overset \sim \to (W,\partial_0 W)
$$
is a simple homotopy equivalence. \item the composite of $h$ with
the inclusion $i\:(W,\partial_0 W) \subset (N,\partial N)$ is
homotopic to $f$. \item The pairs $(W,\partial_1 W)$ and
$(\partial_0 W,\partial_{01}W)$ satisfy the $\pi$-$\pi$ condition.
\end{itemize}
\end{defn}

Observe that the data $(\partial_0 W,h_L)$ define
an embedding up to homotopy in the classical Stallings-Wall sense.

\begin{thm}[Hodgson {\cite[Th.\,\,2.3]{Hodgson}}] \label{Hodgson_embed} Assume $f_K$
is $r$-connected, $r \ge 2k-n+1$, and $k \le n{-}3$. Fix an
embedding up to homotopy
$$
(W_L,h_L)
$$
of the map $f_L\: L \to \partial N$.
Then there exists a relative embedding up to homotopy of $f$
$$
(W,h)
$$
such that $(\partial_0 W,h_L) = (W_L,h_L)$. Furthermore, the
homotopy between $i\circ h$ and $f$ can be taken constant along
$L$.
\end{thm}

\begin{proof}[Proof of \ref{Gamma_connect}:]

By \ref{prototype}, the map $\Gamma$ is surjective on path
components if $3k+2 \le 2n$. We  will show that $\Gamma$ is
surjective on homotopy in degree $j>0$ when $j \le 2n - 3k -2$.
The proof of injectivity is merely a relativized version of the
the proof of surjectivity; to avoid clutter we will omit it.

To prove surjectivity, we can without loss in generality assume
$$
\alpha = E(\beta)\, ,
$$
where $\beta = (M,g)$ is an $n$-thickening.
Hence, $N = M\times [0,1]$.

An element of  $\pi_j(\sec{\bold .}(E_\alpha \to K))$ is represented
by a $j$-parameter family of sections of  the fibration
$$
E_\alpha \to K\,
$$
(by \ref{section_equiv}).
By definition, this is the same thing as a block map
$$
F\:K \times \Delta^j \to \partial N\times \Delta^j
$$
over $\Delta^j$ such that
\begin{itemize}
\item the restriction of $F$ to $K \times \partial\Delta^j$
has the form $F_0 \times \text{id}$;
\item the composition
$$
\begin{CD}
K \times \Delta^j @>F >> \partial N \times \Delta^j @> \subset >>
N \times \Delta^j
\end{CD}
$$
is homotopic over $\Delta^j$ to $f \times \text{id}$, where
the homotopy is held constant along $K \times \partial \Delta^j$.
\end{itemize}
Identify $M$ with the codimension zero
submanifold $M\times 0 \subset \partial N$.
Then we have a codimension zero submanifold
$$
M \times \partial \Delta^j \subset \partial N\times \partial \Delta^j \, ,
$$
a simple homotopy equivalence
$$
\begin{CD}
K \times \partial \Delta^j @> g\times \text{id} > \simeq >
M \times \partial \Delta^j\, ,
\end{CD}
$$
and a map of pairs
$$
\begin{CD}
(K \times \Delta^j,K\times \partial \Delta^j) @> (F,g\times
\text{id}) >>
(\partial N \times \Delta^j,M \times \partial \Delta^j) \\
\end{CD}
$$
By \ref{Hodgson_embed},
there is
\begin{itemize}
\item a codimension zero submanifold
$$
W \subset \partial N \times \Delta^j
$$
whose boundary $\partial W$, when intersected with
 $\partial N \times \partial \Delta^j$,
coincides with $M \times \partial \Delta^j$, and \item a simple
homotopy equivalence
$$
\begin{CD}
K \times \Delta^j @> h > \simeq > W \\
\end{CD}
$$
extending $g\times \text{id}$,
which, when followed by the inclusion into $\partial N \times \Delta^j$
is homotopic to $F$. The homotopy can be assumed
to be fixed along $K \times \partial \Delta^j$.
\end{itemize}

Consider the manifold
$$
V \,\, := \,\, N \times \Delta^j \cup W \times [0,1]        \, ,
$$
where the gluing is along $W \times 0 \subset  \partial N \times \Delta^j$.
This can be thought of in two ways:
\begin{itemize}
\item as diffeomorphic to $N \times \Delta^j$ (since $W \times
[0,1]$ is a collar of $W\times 0$), or \item as an $s$-cobordism
between $W \times 1$ and $\text{closure}(\partial N \times
\Delta^j - W\times 0)$.
\end{itemize}

The  $s$-cobordism theorem gives a diffeomorphism
$$
V \cong W \times [0,1]
$$
relative to $W\times 0$. If we take the mapping cylinder of this
identifications and glue them together, we can regard the result
as an $s$-cobordism $U \subset N\times \Delta^j \times \Delta^1$
between $W \times [0,1]$ and $N \times \Delta^j$. This
$s$-cobordism can be regarded as an element of
$\pi_j(F_n(\alpha))$ which lifts the given homotopy class of the
family of sections. This completes the proof of surjectivity.
\end{proof}

\section{Examples}

The problem with making computations of $\pi_0(T_n(K))$ is that it
generally fails to have the structure of a group. However, we can
obtain {\it upper bounds for the number of elements} in this set
in certain fringe cases.
\bigskip

\subsection*{A lemma about section spaces}
Suppose
$$
K = L \cup_\beta D^k\, .
$$
where $L$ is a based CW complex of dimension $\le k{-}2$, and $k
\ge 2$.

Let
$$
p\:E \to K
$$
be a fibration equipped with section $s\: K \to E$. and suppose
that $p$ is $k$-connected (same $k$ as before). Let $F$ be the
fiber of $p$ at the basepoint.

\begin{lem} With respect to these assumptions
there is an isomorphism
$$
\pi_0(\sec(E\to K)) \cong \pi_k(F) \, .
$$
\end{lem}

\begin{proof} Obstruction theory says the obstructions
to making a section $t\: K \to E$ homotopic to the given section
$s$ live in the cohomology groups
$$
H^*(K;\pi_*(F))\, , \quad  * = 1,2,...
$$
(the coefficients are possibly twisted). By the assumptions, the only
non-trivial group occurs in degree $j=k$, so what this implies is
that we have a bijection
$$
\pi_0(\sec(E\to B)) \cong H^k(K;\pi_k(F))\, .
$$

To compute this, let $C_*(K)$ be the cellular chains on the
universal cover of $K$; this is a complex of free ${\Bbb
Z}[\pi]$-module, where $\pi = \pi_1(K)$. Note that $C_{k-1}(K)$ is
trivial, and $C_k(K) = {\Bbb Z}[\pi]$. The cohomology group is
then easily computed from this and we get $H^k(K;\pi_k(F)) =
\pi_k(F)$.
\end{proof}

\begin{cor}[Stable Version]
With the same assumptions as above we get
$$
\pi_0(\sec^{\rm st}(E\to K)) = \pi_k^{\rm st}(F)
$$
where the right side is the stable homotopy of $F$ in degree $k$.
\end{cor}

\begin{proof} The stable section space in question is the same
thing as the section space of
$$
Q_\bullet E \to K\, .
$$
Apply the the previous lemma and use the fact that $\pi_k^{\rm
st}(F) = \pi_k(QF)$.
\end{proof}

\bigskip

With $K$ as above, we consider $K \times K$. This has the form
$$
K \times K  = (L \times K \cup K\times L) \cup D^{2k}
$$
where the space in parenthesis has dimension $\le 2k - 2$.

Consider the case $2k = n$. By Wall, the tangent map gives an
isomorphism $\pi_0({\bold T}_{n{+}1}(K)) = [K,BO]$.

By ``Corollary C,'' we have a short exact sequence of sets
$$
\pi_1(\sec^{\rm st}({\cal E}_{\xi\oplus \epsilon} \to K \times K))
\to \pi_0({\bold T}_{n}(K)) \overset \tau \to [K,BO] \to 1 \, .
$$
We will compute the first term of this sequence.

The main point is that there is an identification
$$
\Omega \sec^{\rm st}({\cal E}_{\xi\oplus \epsilon} \to K \times K)
\simeq \sec^{\rm st}({\cal E}_{\xi} \to K \times K) \, .
$$

The fiber of the fibration appearing on the right side is of the
form $S^n \smsh (\Omega K)_+$, which is $(2k {-} 1)$-connected.
Hence the projection map for that
 fibration is $(2k)$-connected.

The corollary above then says
$$
\pi_1(\sec^{\rm st}({\cal E}_{\xi\oplus \epsilon} \to K \times K))
\,\, \cong \,\, \pi^{\rm st}_{2k}(S^{n}\smsh (\Omega K)_+)  =
\pi_0^{\rm st}((\Omega K)_+) = {\Bbb Z}[\pi] \, ,
$$
where $\pi = \pi_1(K)$ and the last identification made use of the
Hurewicz theorem.

\begin{cor} With $\pi = \pi_1(K)$, there is a short exact
sequence of sets
$$
{\Bbb Z}[\pi] \to \pi_0({\bold T}_{2k}(K))\overset \tau \to [K,BO]
\to 1 \, .
$$
In particular, $\pi_0({\bold T}_{2k}(K))$ is countable.
\end{cor}

\subsection*{Example 1} In the above,
assume further that $L$ is a CW complex of dimension $\le k{-}2$
with the property that its homology vanishes in all degrees. Then
$K$ has the singular {\it homology of a $k$-sphere} in this case.

Furthermore,
$$
[K,BO] \cong [S^k,BO] \cong \pi_{k{-}1}(O)\, .
$$
(the first isomorphism uses the fact that $BO$ is a loop space and
that $\Sigma K \simeq S^{k{+}1}$).

In particular, if
$$
k \equiv 3,5,6,7  \text{ mod } 8
$$
then $[K,BO] = 0$.

For example, if $k = 11$, then there is a surjective function
$$
\Bbb Z[\pi] \twoheadrightarrow  \pi_0({\bold T}_{22}(K)) \, .
$$

\subsection*{Example 2} In this example
$$
K = S^1 \vee S^{11},
$$

then there is again a short exact sequence of sets
$$
\Bbb Z[t,t^{-1}] \to \pi_0({\bold T}_{22}(S^1 \vee S^{11})) \to
[S^1 \vee S^{11}, BO] \to 1,
$$
and we have
$$
[S^1 \vee S^{11},BO] \cong {\Bbb Z}_2 \, .
$$
Consequently, the sequence of sets
$$
\Bbb Z[t,t^{-1}] \to \pi_0({\bold T}_{22}(S^1 \vee S^{11})) \to
{\Bbb Z}_2 \to 1
$$
is exact.

\end{document}